\documentclass{article}

\usepackage{amsmath,amssymb,amsthm,bm}
\usepackage[a4paper]{geometry}
\usepackage{algorithm,algorithmic}
\usepackage{caption}
\usepackage{graphicx}
\usepackage{hyperref}
\usepackage{xcolor}
\usepackage[T1]{fontenc}
\usepackage{multirow}
\usepackage{nicematrix}
\usepackage{authblk}

\DeclareMathOperator{\Find}{{\tt{find}}}
\DeclareMathOperator{\Length}{{\tt{length}}}
\DeclareMathOperator{\Blkdiag}{{\tt{blkdiag}}}

\newcommand*{\set}[1]{\left\lbrace#1\right\rbrace}

\newcommand*{\herm}{^*}

\newcommand*{\conj}[1]{\overline{#1}}
\newcommand*{\md}{\mathop{}\mathopen{}\mathrm{d}}

\newcommand*{\mi}{\mathrm i}

\def\adots{\mathinner{\mkern2mu\raise1pt\hbox{.}\mkern2mu
    \raise4pt\hbox{.}\mkern2mu\raise7pt\hbox{.}\mkern1mu}}

\newcommand*{\ddiam}{\delta_{\mathrm{eig}}}
\newcommand*{\dres}{\delta_{\mathrm{res}}}
\newcommand*{\dind}{\delta_{\mathrm{ind}}}
\newcommand*{\drex}{\delta_{\mathrm{relax}}}
\newcommand*{\dnr}{\delta_{\mathrm{nr}}}

\newcommand*{\dmax}{d_{\mathrm{max}}}

\newcommand*{\ksub}{k_{\mathrm{sub}}}
\newcommand*{\kin}{k_{\mathrm{in}}}

\newcommand*{\Lambdatmp}{\Lambda_{\mathrm{tmp}}}
\newcommand*{\Vtmp}{V_{\mathrm{tmp}}}

\usepackage[normalem]{ulem}

\graphicspath{{fig/}}

\newenvironment{keywords}{\medskip\textbf{Keywords:}}{}
\newenvironment{AMS}{\medskip\textbf{AMS subject classifications (2020).}}{}

\theoremstyle{plain}

\theoremstyle{definition}

\theoremstyle{remark}
\newtheorem{remark}{Remark}

\definecolor{bkgndcolor}{rgb}{1,1,1}

\title{Solving nonlinear eigenvalue problems via contour integration and region partitioning}
\author[1]{Yuqi Liu}
\author[2]{Jose E. Roman}
\author[3,4]{Meiyue Shao}
\affil[1]{School of Mathematical Sciences, Fudan University, Shanghai 200433,
China}
\affil[2]{D. Sistemes Inform\`atics i Computaci\'o, Universitat Polit\`ecnica de Val\`encia, Cam\'i de Vera s/n, 46022 Val\`encia, Spain}
\affil[3]{School of Data Science, Fudan University, Shanghai 200433, China}
\affil[4]{Shanghai Key Laboratory for Contemporary Applied Mathematics, Fudan
University, Shanghai 200433, China}
\date{\today}

\begin{document}
\pagecolor{bkgndcolor}

\maketitle

\begin{abstract}
In this work, we combine Beyn's method and the recently developed
recursive integral method (RIM) to propose a contour integral-based,
region partitioning eigensolver for nonlinear eigenvalue problems.
A new partitioning criterion is employed to eliminate the need for a
problem-dependent parameter, making our algorithm much more robust compared to
the original RIM.
Moreover, our algorithm can be directly applied to regions containing
singularities or accumulation points, which are typically challenging
for existing nonlinear eigensolvers to handle.
Comprehensive numerical experiments are provided to demonstrate that
the proposed algorithm is particularly well suited
for dealing with regions including many eigenvalues.
\end{abstract}



\section{Introduction}
Given a \(\xi\)-dependent matrix
\(T(\xi)\colon\Omega\rightarrow\mathbb{C}^{n\times n}\), a nonlinear
eigenvalue problem (NEP) is to compute the eigenvector
\(v\in\mathbb{C}^n\setminus\set{0}\) and the corresponding eigenvalue
\(\lambda\in\Omega\) satisfying
\[
T(\lambda)v=0,
\]
where \(\Omega\subseteq\mathbb{C}\) is a connected
region~\cite{GT2017,MV2004}.
In this work, we are particularly interested in the cases where the function
\(T(\xi)\) is \emph{nonlinear} with respect to \(\xi\).
This kind of problems has gained increasing relevance in recent years,
and its applications can be found in various disciplines such as
fluid dynamics~\cite{HR2023,TH2001}, material sciences~\cite{CER2020},
and computational chemistry~\cite{BV2007}.

One of the most prevalent applications of NEP in recent years is the
calculation of quasinormal modes of dispersive materials, which is also
referred to as QNM analysis~\cite{DNG2020,LYG2019,NDZ2023}.
In such applications, magnitudes associated with Maxwell equations are usually discretized by
finite element methods, resulting in NEPs.
The eigenvalues of these NEPs contain the damping and frequency of
the wave's behavior, which is particularly important for understanding
complex materials such as semiconductors and metamaterials.
However, to achieve reasonably good accuracy, researchers usually construct
physical models with high-order polynomials or rational functions,
making numerical computation challenging.

A variety of algorithms can be employed to solve these NEPs,
one class of successful methods is that of the Krylov-based algorithms,
where the NEP is first approximated and reformulated as a linear
eigenvalue problem via companion linearization, which is then solved by
Krylov methods.
This class of methods has undergone significant developments in recent years,
resulting in the emergence of numerous algorithms based on Taylor
expansion~\cite{JMM2012}, Chebyshev interpolation~\cite{KR2014}, and rational
approximation~\cite{CR2021,GVM2014}.

Nevertheless, in many applications, we are concerned with finding all the
eigenvalues within a specific region, as we usually do for QNM.
In such cases, contour integral-based algorithms, which employ the
quadrature rules to approximate the eigenvalues enclosed by a certain
contour, make their way into the spotlight.
Although these algorithms were originally designed for linear
eigenvalue problems~\cite{Polizzi2009,SS2003,SS2007}, some recent developments
have extended them to more general nonlinear
cases~\cite{AST2009,Beyn2012,GMP2018,YS2013}.

However, it is important to note that these algorithms are typically
designed for cases where only a modest number of eigenvalues are to be
computed.
Note that, unlike linear eigenvalue problems, the eigenvectors of NEPs
are not necessarily linearly independent.
As the number of eigenvalues increases,
the matrix consisting of the corresponding eigenvectors is more likely
to become ill-conditioned,
which presents a challenge for both Krylov-based and contour
integral-based algorithms.

To address this challenge, a natural idea is to divide the whole region
into several smaller ones, and deal with them individually to reduce the
number of eigenpairs that need to be handled at each step.
In fact, using the idea of partitioning to solve linear symmetric
eigenvalue problems is already a well-established approach, with robust
mathematical libraries available~\cite{LXE2019}.
This approach employs the theory of spectral density~\cite{LSY2016} to
divide the real axis into several subintervals, and then solve the
eigenvalue problem on each subinterval.

As for the nonlinear case, there are also some algorithms based on a similar
idea.
The reduced subspace iteration (RSI) algorithm~\cite{TS2024} is a recently
developed Krylov-based solver for NEP.
A key distinction between this approach and other Krylov-based algorithms
is that, once it approximates the eigenvalues around a certain shift with
the Arnoldi process, a \(K\)-means algorithm is employed to split the
approximate eigenvalues into multiple clusters.
The cluster cores are then selected as new shifts, and further Arnoldi
processes are implemented on these shifts, respectively, to obtain more
accurate solutions.
Unlike classical Krylov-based nonlinear solvers, appropriate shifts are
selected adaptively without a priori information.
Therefore, RSI is particularly well suited for solving problems involving
many eigenvalues.

In contrast, the recursive integral method (RIM)~\cite{HSS2016} is another
partitioning algorithm, but using contour integration.
This algorithm shares a similar idea to the well-known bisection algorithm
for standard symmetric tridiagonal eigenproblems.
It uses a quadrature rule-based indicator to determine whether there is
any eigenvalue contained in the region of interest.
If so, the region is divided into several subregions.
The same procedure will be repeated recursively until the subregions are
sufficiently small to meet the desired accuracy.
Several enhancements have been made to RIM to improve its practicality and
cost-effectiveness~\cite{HSY2018,HSY2020}.
A nonlinear implementation is also proposed in~\cite{XS2023}.

However, in some practical applications, both RSI and RIM face issues of
missing eigenvalues.
This problem becomes particularly severe when the region of interest contains
a large number of eigenvalues or is close to singularities.
In applications such as analyzing the resonances of a material, missing
eigenvalues can lead to inaccurate predictions.
Therefore, our main goal is to design a more robust algorithm that
overcomes this issue.

In this work, we adopt the idea of the RIM but remove the indicators
used in that method.
Instead, a novel partitioning criterion based on Beyn's method~\cite{Beyn2012}
is employed to estimate the information of the region more accurately.
The proposed algorithm is capable of handling large-scale problems with
many required eigenvalues, robust to singularities, accumulation
points, and high multiplicities, and also less prone to missing eigenvalues.

The remainder of this paper is organized as follows.
In Section~\ref{sec:pre}, we will briefly introduce the two main components
of our algorithm: Beyn's method and the RIM.
Then, in Section~\ref{sec:alg}, we will provide a comprehensive description
of our proposed new algorithm, including implementation details.
Finally, in Section~\ref{sec:numerexp}, we present comprehensive numerical
experiments that demonstrate the excellent performance of our algorithm in
specific physical applications.

For simplicity, we have used MATLAB functions \(\Find\), \(\Length\), and \(\Blkdiag\) in the pseudocode.
Their meanings are as follows:
\begin{itemize}
\item\(\Find\): Returns the indices of the elements in an array that satisfy a given condition.
\item\(\Length\): Returns the length of an array.
\item\(\Blkdiag\): Constructs a block diagonal matrix from the input matrices.
\end{itemize}

\section{Preliminaries}
\label{sec:pre}

\subsection{Beyn's method and quadrature rules}
Beyn's method~\cite{Beyn2012} is one of the most well-known contour
integral-based nonlinear eigensolvers.
To describe this algorithm, we define the zeroth and first moments
\[
\mathcal{M}_0=\frac{1}{2\pi\mi}\int_{\partial\Omega}T(\xi)^{-1}Z\md\xi,\qquad
\mathcal{M}_1=\frac{1}{2\pi\mi}\int_{\partial\Omega}\xi T(\xi)^{-1}Z\md\xi,
\qquad Z\in\mathbb{C}^{n\times k}.
\]
Assuming that there are exactly \(k\) eigenvalues of the NEP
\(T(\lambda)v=0\) lying within the region \(\Omega\), a clever proof provided
in~\cite{Beyn2012} states that these eigenvalues, as well as their
corresponding eigenvectors, can be obtained by performing a spectral
decomposition on the matrix \(V_0\herm\mathcal{M}_1W_0\Sigma_0^{-1}\) for
almost any \(Z\in\mathbb{C}^{n\times k}\), where
\(\mathcal{M}_0=V_0\Sigma_0W_0\herm\) is the singular value decomposition
(SVD) of~\(\mathcal{M}_0\).

However, in numerical implementations, these moments should be approximated
by quadrature rules.
In this work, to facilitate region partitioning, we adopt a rectangular
contour combined with the Gauss--Legendre quadrature rule.
Similarly to the classical approach of using the trapezoidal rule on an
elliptical contour, the Gauss--Legendre rule on each edge of a polygon has
been proven to converge exponentially with respect to the
number of quadrature nodes~\cite{JSS2021,Trefethen2008}.

\begin{figure}[tb!]
\centering
\includegraphics[width=0.8\linewidth]{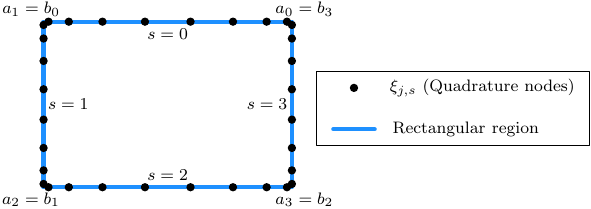}
\caption{Gauss--Legendre quadrature rule on the boundary of a rectangular
region.}
\label{fig:quadraturerule}
\end{figure}

When applying quadrature rules to the boundary of a rectangular region, we
assign \(N/4\) quadrature nodes equally to each of its edges.
Then, the quadrature nodes and the corresponding weights can be obtained
through a transformation of the form
\begin{equation}
\label{eq:gltrans}
\begin{aligned}
&\xi_{j,s}=\frac{\tilde\xi_j\cdot(b_s-a_s)+(b_s+a_s)}{2},
\qquad\omega_{j,s}=\frac{\tilde\omega_j\cdot(b_s-a_s)}{2},\\
&j=0,\dotsc,\frac{N}{4}-1,
\qquad s=0,1,2,3,
\end{aligned}
\end{equation}
where \(\tilde\xi_j\), \(\tilde\omega_j\), \(j=0\), \(\dotsc\), \(N/4-1\) are
the positions and weights of the standard Gauss--Legendre quadrature rule on
the real interval \([-1, 1]\).
(Note that \(\tilde\xi_j\) and \(\tilde\omega_j\) can be calculated with high
accuracy and efficiency by solving a tridiagonal eigenvalue
problem~\cite{Gautschi2012}.)
Also, for convenience, we label the four edges with \(s=0\), \(1\), \(2\),
\(3\) in a counterclockwise order, and \(a_s\), \(b_s\in\mathbb{C}\) are
the starting and ending points of the edge \(s\); see
Figure~\ref{fig:quadraturerule}.

Consequently, when applying to Beyn's method, approximate moments read
\begin{equation}
\label{eq:nummoment}
\mathcal{M}_{0,N}=\frac{1}{2\pi\mi}\sum_{s=0}^{3}
\sum_{j=0}^{N/4-1}\omega_{j,s}T(\xi_{j,s})^{-1}Z,
\qquad\mathcal{M}_{1,N}=\frac{1}{2\pi\mi}\sum_{s=0}^{3}
\sum_{j=0}^{N/4-1}\omega_{j,s}\xi_{j,s}T(\xi_{j,s})^{-1}Z.
\end{equation}
For reference, we list Beyn's method with the Gauss--Legendre
quadrature rule in Algorithm~\ref{alg:beyn}.

\begin{algorithm}[tb!]
\caption{Beyn's method with Gauss--Legendre
quadrature rule on a rectangular contour.}
\begin{algorithmic}[1]
\label{alg:beyn}
\REQUIRE The parameter-dependent matrix
\(T(\xi)\colon\Omega\rightarrow\mathbb{C}^{n\times n}\), the initial
guess \(Z\in\mathbb{C}^{n\times k}\), the number of quadrature nodes
\(N\), and the rectangular region of interest \(\Omega\)
\ENSURE Approximate eigenvalues \(\Lambda\) and eigenvectors \(V\)
\STATE \textbf{function} \([\Lambda,V]={\tt{BeynGL}}(Z,\Omega)\)
\STATE Generate nodes \(\xi_{j,s}\) and weights \(\omega_{j,s}\) on
\(\partial\Omega\) by~\eqref{eq:gltrans} for \(j=0\), \(\dotsc\),
\(N/4-1\), \(s=0\), \(1\), \(2\), \(3\)
\STATE \(\mathcal{M}_{0,N}\gets\frac{1}{2\pi\mi}\sum_{s=0}^{3}
\sum_{j=0}^{N/4-1}\omega_{j,s}T(\xi_{j,s})^{-1}Z\)
\STATE \(\mathcal{M}_{1,N}\gets\frac{1}{2\pi\mi}\sum_{s=0}^{3}
\sum_{j=0}^{N/4-1}\omega_{j,s}\xi_{j,s}T(\xi_{j,s})^{-1}Z\)
\STATE Singular value decomposition \(\mathcal{M}_{0,N}=V_0\Sigma_0 W_0\herm\)
\STATE \(\breve{\mathcal{M}}_{1,N}\gets V_0\herm\mathcal{M}_{1,N}W_0\Sigma_0^{-1}\)
\STATE Eigenvalue decomposition \(\breve{\mathcal{M}}_{1,N}=S\Lambda S^{-1}\)
\STATE \(V\gets V_0S\)
\STATE \textbf{end function}
\end{algorithmic}
\end{algorithm}

\subsection{The recursive integral method}
\label{sub-sec:rim}
The recursive integral method (RIM)~\cite{HSY2018,HSS2016,HSY2020,XS2023} is
a nonlinear eigensolver based on contour integration and region partitioning.
The idea of RIM is that, for a region \(\Omega\subseteq\mathbb{C}\), the
operator
\[
I_{\Omega}(z)=\Bigl\lVert
\frac{1}{2\pi\mi}\int_{\partial\Omega}T(\xi)^{-1}\frac{z}{\lVert z\rVert_2}
\md\xi\Bigr\rVert_2
\]
will return zero for almost any \(z\in\mathbb{C}^n\), \(z\neq0\), if there is
no eigenvalue of \(T(\xi)\) inside \(\Omega\).
Therefore, for a region \(\Omega_*\), we can determine if \(\Omega_*\)
contains an eigenvalue by evaluating \(I_{\Omega_*}(z)\).
The regions with \(I_{\Omega_*}(z)>0\) will be divided recursively until
the diameter of the currently processed region is sufficiently small to meet
the accuracy requirement for the eigenvalues.
At this point, the algorithm returns the center of that region as an
approximate eigenvalue, which is similar to the bisection algorithm
used in standard symmetric tridiagonal eigenproblems.

Of course, the contour integration in \(I_{\Omega}(z)\) is by no means
calculated mathematically.
RIM employs an indicator with an \(N\)-point quadrature rule
\begin{equation}
\label{eq:rimind}
I_{\Omega,N}(z)=\Bigl\lVert\frac{1}{2\pi\mi}\sum_{j=0}^{N-1}\omega_j
T(\xi_j)^{-1}\frac{z}{\lVert z\rVert_2}\Bigr\rVert_2,
\end{equation}
where \(\xi_j\)'s and \(\omega_j\)'s are some quadrature nodes and weights
on \(\partial\Omega\).
A classification threshold \(\dind>0\) is used to determine partitioning---whenever
\(I_{\Omega_*,N}(z)\ge\dind\), we suppose that \(\Omega_*\)
contains at least one eigenvalue, and divide it into several subregions.
The framework of the RIM is listed in Algorithm~\ref{alg:rim}.

\begin{algorithm}[tb!]
\caption{Recursive integral method (RIM)}
\begin{algorithmic}[1]
\label{alg:rim}
\REQUIRE Rectangular region of interest \(\Omega\), accuracy tolerance for
eigenvalues \(\ddiam\), classification threshold for the indicator \(\dind\),
and the number of quadrature nodes \(N\)
\ENSURE Approximate eigenvalues \(\Lambda\)
\STATE \textbf{function} \([\Lambda]={\tt{RIM}}(\Omega)\)
\STATE \(\Lambda=[~]\)
\STATE Randomly generate \(z\in\mathbb{C}^n\) with \(\lVert z\rVert_2=1\)
\STATE Generate nodes \(\xi_{j,s}\) and weights \(\omega_{j,s}\) on
\(\partial\Omega\) by~\eqref{eq:gltrans} for \(j=0\), \(\dotsc\),
\(N/4-1\), \(s=0\), \(1\), \(2\), \(3\)
\STATE Compute \(I_{\Omega,N}(z)=\bigl\lVert\sum_{s=0}^{3}
\sum_{j=0}^{N/4-1}\omega_{j,s}T(\xi_{j,s})^{-1}z\bigr\rVert_2/(2\pi)\)
\label{alg-step:rim-Iz}
\IF{\(I_{\Omega,N}(z)<\dind\)}
  \STATE \textbf{exit} \quad\COMMENT{No eigenvalues in \(\Omega\), return void}
\ELSIF{the diameter of \(\Omega\) is larger than \(2\cdot\ddiam\)}
  \STATE Divide \(\Omega\) into \(\Omega_j\) for \(j=0\), \(\dotsc\), \(3\)
  \quad\COMMENT{by Figure~\ref{fig:dividepattern}}
  \FOR{\(j=0\), \(\dotsc\), \(3\)}
    \STATE \(\Lambda\gets\Blkdiag(\Lambda,{\tt{RIM}}(\Omega_j))\)
    \quad\COMMENT{Not converged yet, divide the region}
  \ENDFOR
\ELSE
  \STATE \(\Lambda\gets\) the center of \(\Omega\)
  \quad\COMMENT{Converged, return center of \(\Omega\) as approximate eigenvalue}
  \label{alg-step:converge}
\ENDIF
\STATE \textbf{end function}
\end{algorithmic}
\end{algorithm}

An advantage of this algorithm is that RIM is able to isolate all the
eigenvalues without a priori spectral information~\cite{HSS2016}.
The algorithm can locate the eigenvalues even with only four quadrature
nodes, which makes it a powerful tool to solve transmission eigenvalue
problems~\cite{GST2022,HSS2016}.

However, the classification threshold for the indicator \(\dind\) used in Algorithm~\ref{alg:rim}
is highly dependent on the specific problem, making it a tricky task to
choose \(\dind\) in practice~\cite{XS2023}.
A too large \(\dind\) may result in the loss of certain eigenvalues, while
a too small one could lead to unnecessary computation of regions that do
not actually contain an eigenvalue.
Even for the same problem, as subregions become finer during
partitioning, the return value of \(I_{\Omega,N}(z)\) may vary significantly,
making the parameter ineffective and the algorithm inefficient. 

Furthermore, the indicator~\eqref{eq:rimind} operates on a single vector
\(z\), and therefore can only return binary information about the existence or not
of eigenvalues.
Thus, the algorithm has to increase the search depth in order to ensure that
only one eigenvalue is contained in each subregion in the end, which can
exponentially increase the computational cost.

\begin{remark}[How to divide the region]
In this paper, we always assume that the region of interest is rectangular.
When partitioning, the region is divided into four equal smaller
rectangular subregions in a cross-shaped manner; see
Figure~\ref{fig:dividepattern}.
For convenience, we label these subregions by \(0\), \(1\), \(2\), \(3\).
\begin{figure}[tb!]
\centering
\includegraphics[width=0.6\linewidth]{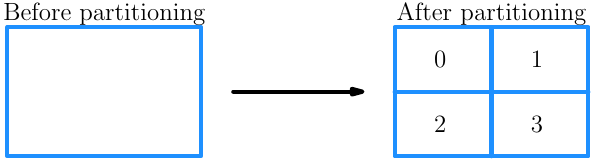}
\caption{Partitioning pattern.}
\label{fig:dividepattern}
\end{figure}
\end{remark}

\section{Proposed algorithm}
\label{sec:alg}

To improve the RIM, we have two primary goals:
to avoid excessive search depth and to eliminate the problem-dependent
threshold \(\dind\).

\subsection{RIM by numerical rank}
To achieve the first goal, we would like to employ a matrix
\(Z\in\mathbb{C}^{n\times\ksub}\) for some \(\ksub>1\), instead of
using a single vector \(z\in\mathbb{C}^n\) in~\eqref{eq:rimind}.
Note that it has been proved in~\cite{Beyn2012} that for almost any \(Z\) the
rank of the zeroth moment
\[
\mathcal{M}_0=\frac{1}{2\pi\mi}\int_{\partial\Omega}T(\xi)^{-1}Z\md\xi
\]
is exactly \(\min\{\ksub,k\}\), where \(k\) is the number of eigenvalues of
\(T(\xi)\) within \(\Omega\).
Thus, by estimating the rank of \(\mathcal{M}_0\), we can determine whether
there are eigenvalues within the region and also estimate the number of them.

Therefore, we can divide the region only if \(\mathcal{M}_0\) is of full rank.
When \(\mathcal{M}_0\) is rank deficient, it indicates that the region
contains less than \(\ksub\) eigenvalues, and we can apply existing NEP
solvers to compute them directly to reduce the excessive search depth which
may emerge in RIM.%
\footnote{A similar idea can be found in~\cite{XS2023}, where Beyn's method
is employed to terminate the RIM earlier.}
This allows the algorithm to stop earlier, without the need to divide the
region into excessively small pieces.

To estimate the rank of \(\mathcal{M}_0\) numerically, a typical approach is
to perform a truncated SVD on the approximate moment~\(\mathcal{M}_{0,N}\)
(see~\eqref{eq:nummoment}) by dropping singular values below a classification
threshold~\(\dnr\).
The rank estimated by this truncated SVD procedure is then taken as the
numerical rank of \(\mathcal{M}_0\).
We summarize this modification of RIM in Algorithm~\ref{alg:rim_numrank}.

\begin{algorithm}[tb!]
\caption{The recursive integral method by numerical rank}
\begin{algorithmic}[1]
\label{alg:rim_numrank}
\REQUIRE Rectangular region of interest \(\Omega\), classification threshold for numerical rank \(\dnr\), relaxation parameter
\(\drex\), maximum number of eigenvalues of each subregion \(\ksub\), the
number of quadrature nodes \(N\), and a high performance NEP solver
\({\tt{NEPsolver}}()\) for small-scale problems
\ENSURE Approximate eigenvalues \(\Lambda\)
\STATE \textbf{function} \([\Lambda]={\tt{RIMNR}}(\Omega)\)
\STATE \(\Lambda=[~]\)
\STATE Randomly generate \(Z\in\mathbb{C}^{n\times\ksub}\)
\STATE Generate nodes \(\xi_{j,s}\) and weights \(\omega_{j,s}\) on
\(\partial\Omega\) by~\eqref{eq:gltrans} for \(j=0\), \(\dotsc\),
\(N/4-1\), \(s=0\), \(1\), \(2\), \(3\)
\STATE \(\mathcal{M}_{0,N}\gets\frac{1}{2\pi\mi}\sum_{s=0}^{3}
\sum_{j=0}^{N/4-1}\omega_{j,s}T(\xi_{j,s})^{-1}Z\)
\STATE Singular value decomposition \(\mathcal{M}_{0,N}=V_0\Sigma_0W_0\herm\)
\label{alg-step:rim_numrank_svd}
\STATE \(\kin\gets\Length\bigl(\Find({\tt{diag}}(\Sigma_0)>\dnr)\bigr)\)
\quad\COMMENT{The number of singular values larger than \(\dnr\)}
\IF{\(\kin/\ksub<\drex\)}
  \STATE \(\Lambda\gets{\tt{NEPsolver}}(\Omega)\)
  \quad\COMMENT{Few eigenvalues in \(\Omega\), use {\tt{NEPsolver}}}
\ELSE
  \STATE Divide \(\Omega\) into \(\Omega_j\) for \(j=0\), \(\dotsc\), \(3\)
  \quad\COMMENT{by Figure~\ref{fig:dividepattern}}
  \FOR{\(j=0\), \(\dotsc\), \(3\)}
    \STATE \(\Lambda\gets\Blkdiag(\Lambda,{\tt{RIMNR}}(\Omega_j))\)
    \quad\COMMENT{Too many eigenvalues, divide the region}
  \ENDFOR
\ENDIF
\STATE \textbf{end function}
\end{algorithmic}
\end{algorithm}

Nevertheless, this approach shares a similar problem with RIM:
For different problems, the proper threshold \(\dnr\) can be very different,
and a poorly chosen \(\dnr\) leads to a completely wrong estimate of the rank
of \(\mathcal{M}_0\).

\subsection{Recursive Beyn's method}
To address the issue of RIM, we have to introduce a more quantifiable
criterion to determine whether a region should be further divided.
Beginning from Line~\ref{alg-step:rim_numrank_svd} of
Algorithm~\ref{alg:rim_numrank}, this time we continue on Beyn's
method, i.e., Algorithm~\ref{alg:beyn}, until the approximate
eigenvalues \(\Lambda\) are obtained.
Then, instead of using the numerical rank, we employ the number of
\(\Lambda\) lying within the region~\(\Omega\), as the approximate rank of
\(\mathcal{M}_0\).
The new algorithm, referred to as recursive Beyn's method, is listed in
Algorithm~\ref{alg:rbm}.

This approach eliminates the necessity for a classification threshold such as
\(\dind\) and \(\dnr\), which are highly problem-dependent and not robust,
resulting in a more flexible partitioning criterion.
Furthermore, it can be observed from Algorithm~\ref{alg:beyn} that, by the
time \(\mathcal{M}_{0,N}\) is computed, the majority of the computational cost
of Beyn's method has already been incurred.
The remaining steps to complete a full run of Beyn's method are of low computational
complexity, which means the overhead introduced is very low.

\begin{algorithm}[tb!]
\caption{The recursive Beyn's method}
\begin{algorithmic}[1]
\label{alg:rbm}
\REQUIRE Rectangular region of interest \(\Omega\), accuracy tolerance for eigenpairs
\(\dres\), maximum number of eigenvalues of each subregion \(\ksub\),
relaxation parameter \(\drex\), the number of quadrature nodes \(N\),
and the maximum search depth \(\dmax\)
\ENSURE Approximate eigenvalues \(\Lambda\) and eigenvectors \(V\)
\STATE \textbf{function} \([\Lambda,V]={\tt{RBM}}(\dmax,\Omega)\)
\STATE \(V\gets[~]\), \(\Lambda\gets[~]\)
\STATE Randomly generate \(Z\in\mathbb{C}^{n\times \ksub}\)
\STATE \([\Lambdatmp,\Vtmp]={\tt{BeynGL}}(Z,\Omega)\)
\COMMENT{{\tt{BeynGL}} is defined in Algorithm~\ref{alg:beyn}}
\label{alg-step:rbmbeyn}
\STATE \({\tt{idx}}=\Find(\Lambdatmp\in\Omega)\)
\quad\COMMENT{Identify eigenpairs in \(\Omega\)}
\STATE\(\Lambdatmp\gets\Lambdatmp({\tt{idx}},{\tt{idx}})\),
\(\Vtmp\gets\Vtmp(~:~,{\tt{idx}})\)
\STATE \(\kin\gets\Length({\tt{idx}})\)
\STATE Check convergence on \([\Lambdatmp,\Vtmp]\)
\IF{\(\kin/\ksub<\drex\) \AND \([\Lambdatmp,\Vtmp]\) converged}
\label{alg-step:ddiv}
  \STATE \(V\gets\Vtmp\), \(\Lambda\gets\Lambdatmp\)
  \quad\COMMENT{Converged, save the solutions}
\ELSIF{\(\dmax>0\)}
  \STATE Divide \(\Omega\) into \(\Omega_j\) for \(j=0\), \(\dotsc\), \(3\)
  \quad\COMMENT{by Figure~\ref{fig:dividepattern}}
  \FOR{\(j=0\), \(\dotsc\), \(3\)}
    \STATE \([\Lambdatmp,\Vtmp]={\tt{RBM}}(\dmax-1,\Omega_j)\)
    \STATE \(V\gets[V,\Vtmp]\), \(\Lambda\gets\Blkdiag(\Lambda,\Lambdatmp)\)
    \COMMENT{Not converged, divide the region}
  \ENDFOR
\ENDIF
\STATE \textbf{end function}
\end{algorithmic}
\end{algorithm}

\begin{remark}
In both Algorithms~\ref{alg:rim_numrank} and~\ref{alg:rbm}, a
relaxation parameter \(0<\drex\le1\) is introduced.
We use \(\kin/\ksub\ge\drex\) instead of \(\kin=\ksub\) as the criterion to
decide whether to proceed with an additional level of partitioning.
A \(\drex<1\) can help avoid the potential risk of missing eigenvalues in some
challenging cases, at the cost of increasing the search depth of the
algorithm.
Note that, unlike \(\dind\) in RIM, the value used for the parameter \(\drex\)
will not influence the correctness of the algorithm.
\end{remark}

\begin{remark}
In Algorithm~\ref{alg:rbm}, we introduce another parameter \(\dmax\) to limit
the depth of recursion.
The algorithm stops when the maximum search depth \(\dmax\) is reached.
Unexploited subregions can be marked for future processing, if necessary.
This enables the algorithm to deal with difficult regions such as those
containing singularities and eigenvalues with high multiplicities.
\end{remark}

\begin{remark}
It is worth noting that Algorithm~\ref{alg:rim}, as originally defined
in~\cite{HSS2016}, as well as Algorithm~\ref{alg:rim_numrank},
can only return eigenvalues.
Further computations should be done if the eigenvectors are also required.
In contrast, Algorithm~\ref{alg:rbm} can return the corresponding eigenvectors at the same time, since Beyn's method is used.
\end{remark}

\subsection{Accelerate the linear solver with infinite GMRES}

Infinite GMRES (infGMRES)~\cite{JC2022} is a recently proposed algorithm
for solving parameterized linear systems
\[
T(\xi)x=b,\qquad
b\in\mathbb{C}^{n},\qquad
T(\xi)\colon\Omega\rightarrow\mathbb{C}^{n\times n},\qquad
\Omega\subseteq\mathbb{C},
\]
on \(\xi=\xi_0\), \(\xi_1\), \(\dotsc\), altogether efficiently.
In brief, infGMRES uses Taylor expansion and companion linearization to
transform the problem into solving another parameterized linear system
\[
(I-\xi\mathcal{L}_0^{-1}\mathcal{L}_1)\tilde y=\tilde b,\qquad
\tilde b\in\mathbb{C}^{mn},\qquad
\mathcal{L}_0,\mathcal{L}_1\in\mathbb{C}^{mn\times mn},
\]
on \(\xi=\xi_0\), \(\xi_1\), \(\dotsc\), where \(m\in\mathbb{N}\) is a
parameter related to the order of the Taylor expansion.
Since the Krylov subspaces generated by
\(I-\xi_j\mathcal{L}_0^{-1}\mathcal{L}_1\) with different \(\xi_j\)'s are
the same, the linear systems can be solved by a single Arnoldi process.
In~\cite{LRS2024}, we have employed infGMRES to accelerate the contour
integral-based nonlinear eigensolvers.
This idea can also be applied here to enhance the efficiency of
recursive Beyn's method.

However, it is proved in~\cite{JC2022} that the accuracy of infGMRES is
influenced by the distance between the points to be solved and the point
employed for the Taylor expansion (that is, the expansion point).
In order to ensure accuracy, it is often necessary to distribute multiple
expansion points.
To this end, the simplest idea is to first place a uniform grid over the
region \(\Omega\), use these grid points as expansion points to employ
infGMRES, and store the corresponding Arnoldi and Hessenberg matrices.
Then, whenever a linear system is to be solved in Line~\ref{alg-step:rbmbeyn},
Algorithm~\ref{alg:rbm}, we simply choose the nearest expansion point, and
utilize the precomputed infGMRES results at that point for efficient solving.
In Section~\ref{sub-sec:numexp_infgmres}, several numerical examples are
provided to illustrate the performance of employing infGMRES in this way.

We remark that our application of infGMRES here is just heuristic.
In addition to the distance we have just mentioned, the accuracy of infGMRES
is also influenced by multiple factors, such as the singularities of \(T\).
As a result, in practical applications, it is often difficult to determine in
advance how to distribute the expansion points to ensure that all linear
systems are solved accurately.

\section{Numerical experiments}
\label{sec:numerexp}

In this section, we employ recursive Beyn's method
(Algorithm~\ref{alg:rbm}) to solve several NEPs in physical applications.
All numerical experiments were carried out using MATLAB R2023b on a Linux
server with two 16-core Intel Xeon Gold 6226R 2.90 GHz CPUs and 1024~GB of
main memory.
The convergence criterion for an approximate nonlinear eigenpair
\((\hat\lambda,\hat v)\) is
\[
\lVert T(\hat\lambda)\hat v\rVert_2\le
\dres\cdot\lVert T(\hat\lambda)\rVert_2\lVert\hat v\rVert_2,
\]
where \(\dres\) is a user-specified tolerance.
Unless otherwise stated, all test results presented in this section have
achieved an accuracy of \(\dres=10^{-12}\), and the relaxation parameter
\(\drex\) in Algorithm~\ref{alg:rbm} is always set to \(0.8\).

\subsection{Test problems}
\label{sub-sec:overall}
We begin with a brief introduction to three test NEPs, describe their
respective properties and challenges, and present the overall performance of
our algorithm on these problems.

\subsubsection{The {\tt{gun}} problem}
The {\tt{gun}} problem~\cite{LBL2010} is one of the most popular test problems
for nonlinear eigensolvers~\cite{TS2024,VanBeeumen2015}.
The problem is derived from the electromagnetic modeling of waveguide
loaded accelerator cavities.
With a nonlinear boundary condition, the physical problem can be discretized
by the finite element method to
\begin{equation}
\label{eq:gun}
\begin{aligned}
&\Bigl(A_0-\lambda A_1+\mi\sqrt{\lambda-\kappa_{c,1}^2}A_2
+\mi\sqrt{\lambda-\kappa_{c,2}^2}A_3\Bigr)v=0,\\
&\kappa_{c,1},\kappa_{c,2}\in\mathbb{R},
\qquad A_0,A_1,A_2,A_3\in\mathbb{R}^{9956\times 9956}.
\end{aligned}
\end{equation}
In this work, we generate the {\tt{gun}} problem from the NLEVP
collection~\cite{BHM2013}, and take the region of interest to be the square
centered at \((250^2,0)\) with side length
\(2\cdot(300^2-200^2)\)~\cite{VanBeeumen2015}.

\begin{figure}[tb!]
\centering
\includegraphics[width=0.7\linewidth]{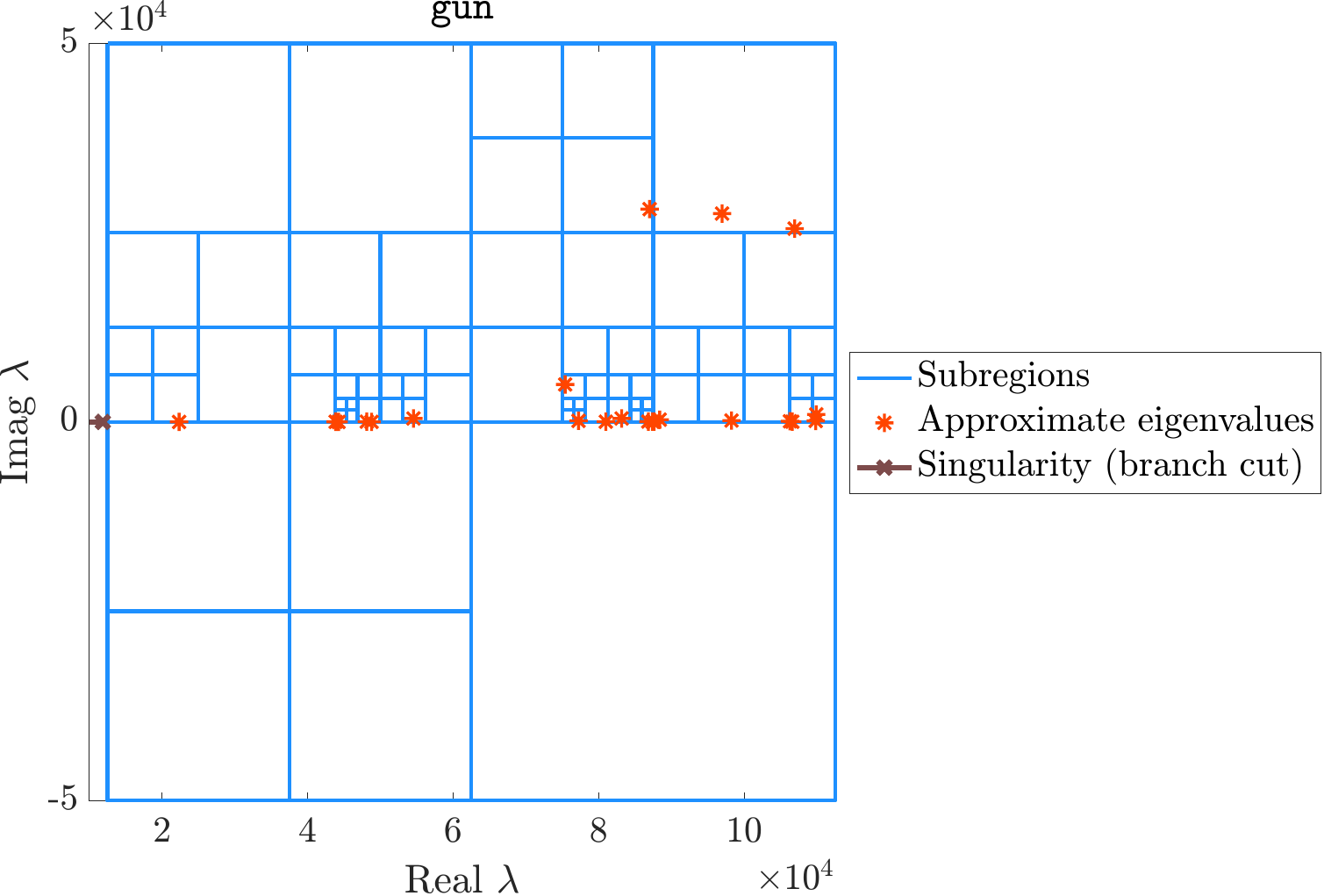}
\caption{Using recursive Beyn's method to solve the {\tt{gun}} problem.
The region of interest is
\([1.25\times10^4,1.125\times10^5]\times[-5\times10^4,5\times10^4]\)
in the complex plane.
The parameters are set as: maximum eigenvalues in each subregion \(\ksub=5\),
number of quadrature nodes on each segment \(N=32\), and maximum search depth
\(\dmax=6\).
All \(22\) eigenvalues as well as their corresponding eigenvectors are found.
The branch cut \((-\infty,108.8774^2]\) is indicated by a line with a
cross mark on its endpoint.}
\label{fig:gun}
\end{figure}

The result of solving the problem using recursive Beyn's method can be found
in Figure~\ref{fig:gun}.
The grid of subregions and the approximate eigenvalues are shown in the figure.
All \(22\) eigenpairs in the region of interest are found and computed with a
relative residual below \(10^{-12}\).
Note that the square roots in~\eqref{eq:gun} are both defined by the principal
branch and therefore lead to a branch cut
\((-\infty,\kappa_{c,2}^2]\) lying close to the boundary of the region.
This introduced some challenges to the subregions on the left-hand side.
However, our algorithm successfully identified all the eigenvalues
accurately.

\subsubsection{The {\tt{photonics\_1}} problem}
\begin{figure}[tb!]
\centering
\includegraphics[width=0.6\linewidth]{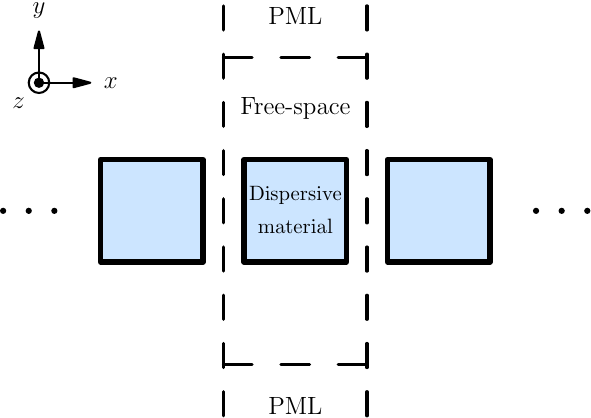}
\caption{Geometric distribution of the {\tt{photonics\_1}} and
{\tt{photonics\_2}} problems.}
\label{fig:geodrude}
\end{figure}

To better demonstrate the performance of our algorithm in practical
applications, we provide two recently proposed problems from quasinormal modes
analysis (QNM).
The physical model~\cite{DNG2020,GDJ2017} being considered is illustrated in
Figure~\ref{fig:geodrude}, an infinite number of rods
are evenly distributed along the \(x\)-axis, and their gaps are filled by
free-space.
Then, along the \(y\)-axis, the unbounded domain is handled using the standard
Cartesian Perfectly Matched Layer (PML) technique.
Additionally, for simplicity, the system is assumed to be invariant along the
\(z\)-axis.

\begin{figure}[tb!]
\centering
\includegraphics[width=0.7\linewidth]{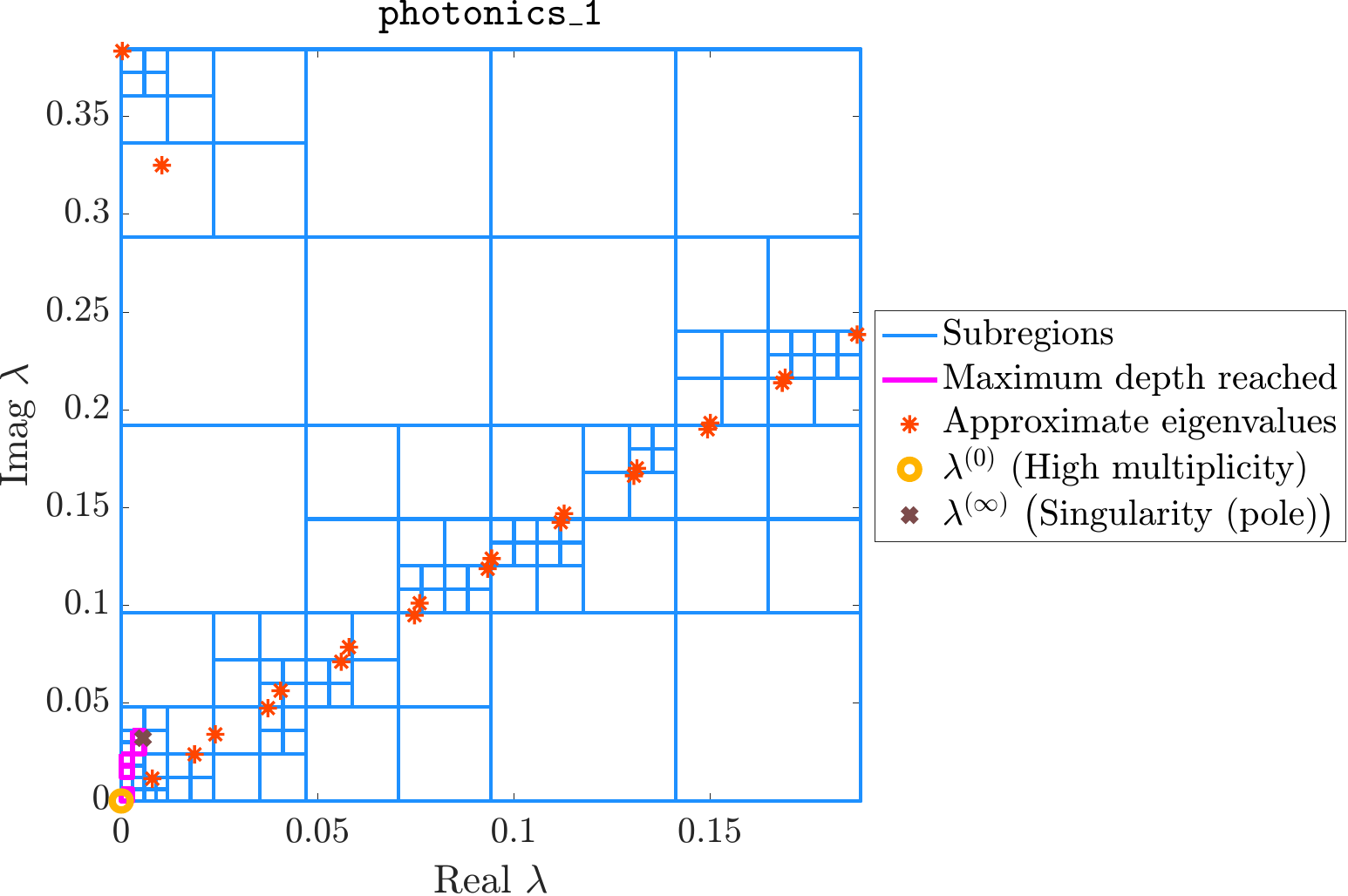}
\caption{Using recursive Beyn's method to solve the {\tt{photonics\_1}}
problem.
The region of interest is \([0,0.188365]\times[0,0.384419]\) in the complex plane.
The parameters are set as: maximum eigenvalues in each subregion \(\ksub=5\),
number of quadrature nodes on each segment \(N=32\), and maximum search depth
\(\dmax=6\).
A total of \(22\) different nonzero eigenvalues as well as their
corresponding eigenvectors are found.
The eigenvalue with high multiplicity \(\lambda^{(0)}\) and the pole
\(\lambda^{(\infty)}\) are also marked in the figure.}
\label{fig:photonics1}
\end{figure}

For the {\tt{photonics\_1}} problem, we take the dispersive material to be
gold, and employ the permittivity model from~\cite{GDJ2017}.
Then, the expression for the permittivity reads
\[
\varepsilon_{\rm{p}}(\mi\omega)=\varepsilon_{\infty}
+\frac{\alpha_1}{(\mi\omega)-\beta_1}
+\frac{\conj{\alpha}_1}{(\mi\omega)-\conj{\beta}_1},
\qquad\varepsilon_{\infty},\alpha_1,\beta_1\in\mathbb{C},
\]
where \(\omega\) is the (complex) frequency of the electric wave.
Using the finite element method, we can discretize Maxwell's equations
on the system to obtain the rational eigenvalue problem
\begin{equation}
\label{eq:polenep}
\bigl(A_0+\lambda^2A_1+\lambda^2\cdot\varepsilon_{\mathrm{p}}(\lambda)A_2\bigr)v=0,
\qquad A_0,A_1,A_2\in\mathbb{C}^{20363\times 20363}.
\end{equation}
Solving~\eqref{eq:polenep} by recursive Beyn's method on the region of
interest \([0,0.188365]\times[0,0.384419]\) leads to
Figure~\ref{fig:photonics1}.

There are two main challenges to this problem:
\begin{enumerate}
\item \(\varepsilon_{\rm{p}}(\lambda^{(0)})=0\) or \(\lambda^{(0)}=0\):
This means that the permittivity of the material or the electric field itself
becomes zero when \(\lambda=\lambda^{(0)}\).
Physically, this solution is trivial.
However, numerically, this eigenvalue usually has an extremely high algebraic
multiplicity, which can cause common numerical algorithms to waste
excessive computational resources at this meaningless point~\cite{DNG2020}.
\item \(\varepsilon_{\rm{p}}(\lambda^{(\infty)})=\infty\):
This is a pole of the rational function \(\varepsilon_{\rm{p}}\).
Numerical solutions may be inaccurate near this point.
\end{enumerate}

In addition to the difficulties mentioned above, there is another
common challenge for various contour integral-based eigensolvers ---
the eigenvalues near the contour.
When such eigenvalues exist, either inside or outside the region, the
quadrature rules are typically observed to be inaccurate~\cite{JSS2021},
resulting in the risk of missing eigenvalues or computing eigenvalues with
large errors/slow convergence.

Note from Figure~\ref{fig:photonics1} that five subregions are marked as
unexploited when the maximum depth~\(6\) is reached (marked in magenta).
Although three of them are due to numerical errors, the other two subregions
capture the two difficult points, \(\lambda^{(0)}\) and
\(\lambda^{(\infty)}\), respectively.
However, we remark that the eigenvalues lying close to \(\lambda^{(0)}\),
\(\lambda^{(\infty)}\), or the boundary of the region of interest,
are still calculated with a relative residual below \(10^{-12}\).

\subsubsection{The {\tt{photonics\_2}} problem}
For the {\tt{photonics\_2}} problem, we consider an even more complicated
model from~\cite{DNG2020}.
With the same geometric model as in Figure~\ref{fig:geodrude}, we use a Drude model for the permittivity
--- a real-coefficient
rational function with respect to \(\mi\omega\) as
\begin{equation}
\label{eq:drudeepsilon}
\varepsilon_{\rm{d}}(\mi\omega)=\varepsilon_{\infty}
-\frac{\omega_{\rm{d}}^2}{-(\mi\omega)^2+\gamma_{\rm{d}}(\mi\omega)},
\qquad\varepsilon_{\infty},\omega_{\rm{d}},\gamma_{\rm{d}}\in\mathbb{R},
\end{equation}
where \(\omega\) is the (complex) frequency of the electric wave.

\begin{figure}[tb!]
\centering
\includegraphics[width=0.8\linewidth]{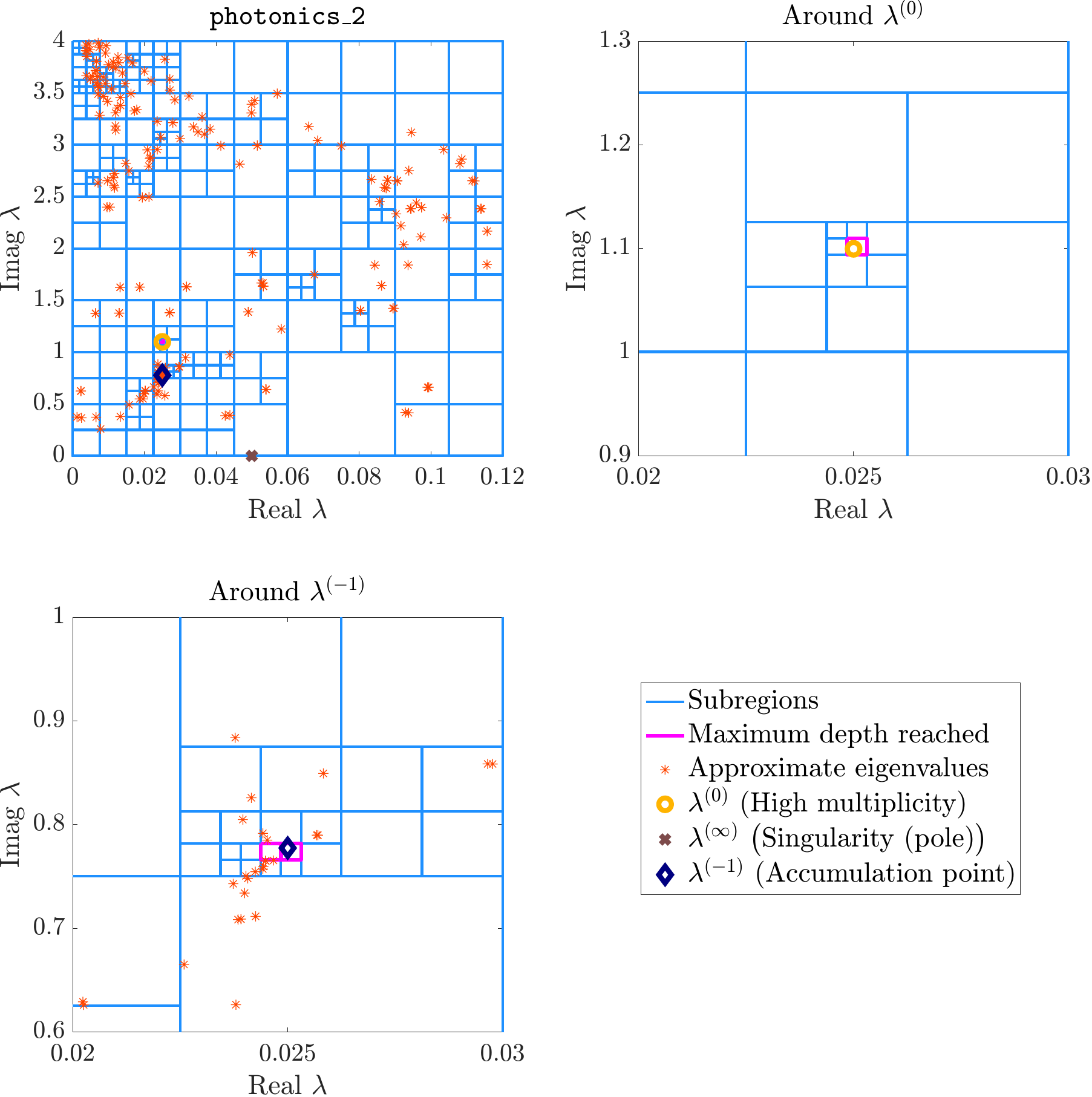}
\caption{Using recursive Beyn's method to solve the {\tt{photonics\_2}}
problem.
The region of interest is \([0,0.12]\times[0.0005,4]\) in the complex plane.
The parameters are set as: maximum eigenvalues in each subregion \(\ksub=20\),
number of quadrature nodes on each segment \(N=512\), and maximum search depth
\(\dmax=8\).
A total of  \(197\) different eigenvalues as well as their corresponding eigenvectors
are found.
The eigenvalue with high multiplicity \(\lambda^{(0)}\), the plasmonic
resonance \(\lambda^{(-1)}\), and the pole \(\lambda^{(\infty)}\) are also
marked in the figure.
To provide a clear picture, we zoom in near \(\lambda^{(0)}\) and
\(\lambda^{(-1)}\).}
\label{fig:demesy5}
\end{figure}

Employing all the parameters used in~\cite{DNG2020} with a finite-element mesh
with \(1183\) elements, we obtain the NEP
\begin{equation}
\label{eq:drudenep}
\bigl(A_0+\lambda^2A_1+\lambda^2\cdot\varepsilon_{\mathrm{d}}(\lambda)A_2\bigr)v=0,
\qquad A_0,A_1,A_2\in\mathbb{C}^{6667\times 6667}.
\end{equation}
This time, to analyze the quasinormal modes of this system, we aim
to compute the eigenvalues within \(\lambda\in[0,0.1]\times[0.0005,4]\)
and their corresponding eigenvectors \(v\).
The output of our algorithm on this problem is illustrated in
Figure~\ref{fig:demesy5}.

The problem {\tt{photonics\_2}} is more difficult than {\tt{photonics\_1}},
although the dimension of the former is only about a third of the latter.
This is because the region of interest now contains almost \(200\)
different eigenvalues.
And, in addition to the previously mentioned difficulties with
\(\lambda^{(0)}\) (where \(\varepsilon_{\rm{d}}(\lambda^{(0)})=0\) or
\(\lambda^{(0)}=0\)) and \(\lambda^{(\infty)}\) (where \(\varepsilon_{\rm{d}}
(\lambda^{(\infty)})=\infty\)),
there is also:
\begin{enumerate}
\setcounter{enumi}{2}
\item \(\varepsilon_{\rm{d}}(\lambda^{(-1)})=-1\):
\(\lambda^{(-1)}\) is an accumulation point of the NEP~\eqref{eq:drudenep}.
The eigenvalues concentrated around it correspond to the plasmonic
resonances, which means their eigenvectors, or electric field, show a
significant oscillation on the interface of the Drude material.
These eigenpairs are also not of interest.
\end{enumerate}

It can be observed from the figure that, with a maximum search depth of \(8\),
both \(\lambda^{(0)}\) and \(\lambda^{(-1)}\) are located by the algorithm.
Although the pole \(\lambda^{(\infty)}=0.05\) lies close to the boundary of
the region of interest, \(197\) eigenpairs are computed with a relative
residual below \(10^{-12}\).

\subsection{Comparison among different partitioning criteria}
\begin{figure}[tb!]
\centering
\includegraphics[width=0.5\linewidth]{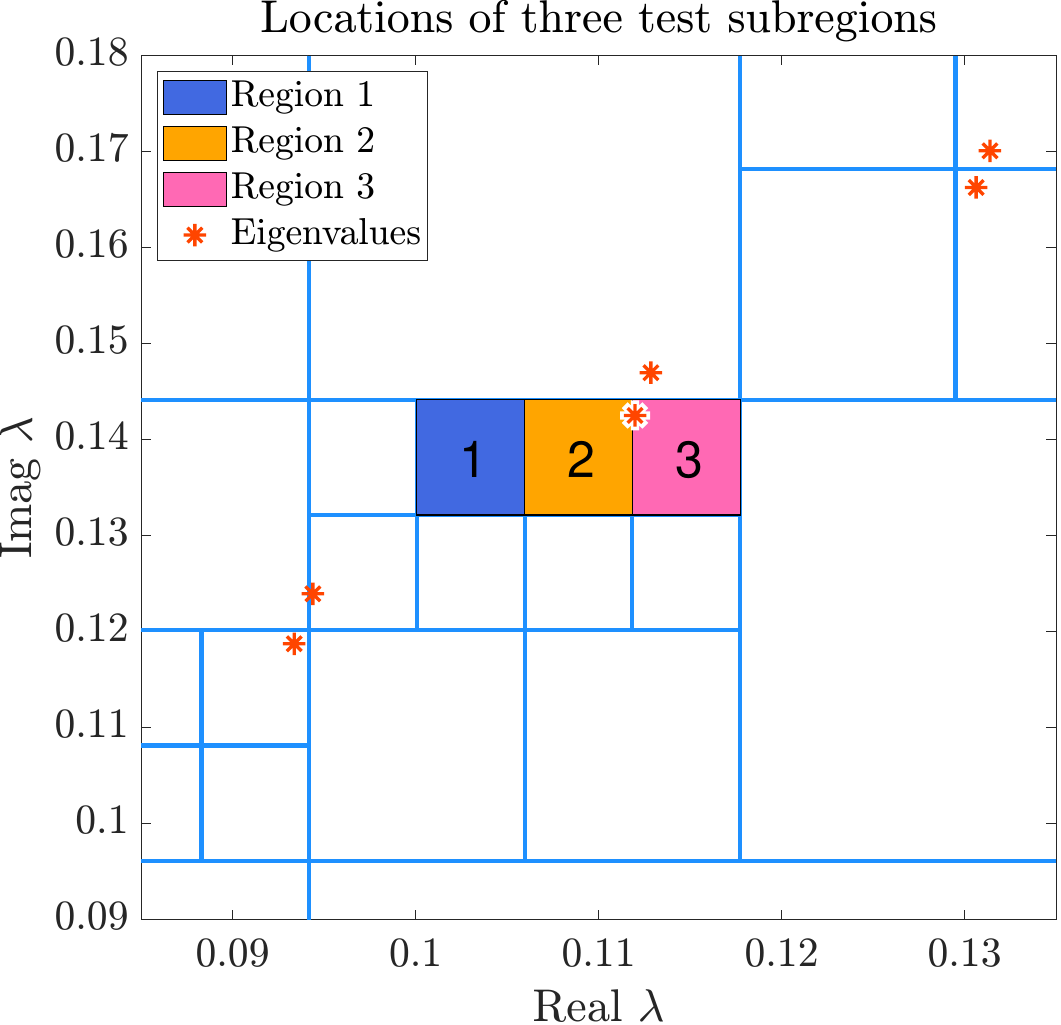}
\caption{Three test regions.
All these regions are chosen from subregions in Figure~\ref{fig:photonics1}.}
\label{fig:3indicator_rg}
\end{figure}

In this subsection, we test whether the proposed algorithm can accurately
reflect the number of eigenvalues within a region.
For comparison, recursive Beyn's method is applied to three different
test regions on the {\tt{photonics\_1}} problem; see
Figure~\ref{fig:3indicator_rg}.
The three test regions are all taken from the subregions generated during the
process of recursive Beyn's method, i.e., from
Figure~\ref{fig:photonics1}.

There are several differences across these regions.
Region \(3\) contains an eigenvalue, while Regions \(1\) and \(2\) do
not contain any eigenvalue.
Another important difference is that there is an eigenvalue lying
extremely close to the boundary between Regions \(2\) and \(3\), posing
a challenge for the contour integration on Regions \(2\) and \(3\).

\begin{figure}[tb!]
\centering
\includegraphics[width=0.8\linewidth]{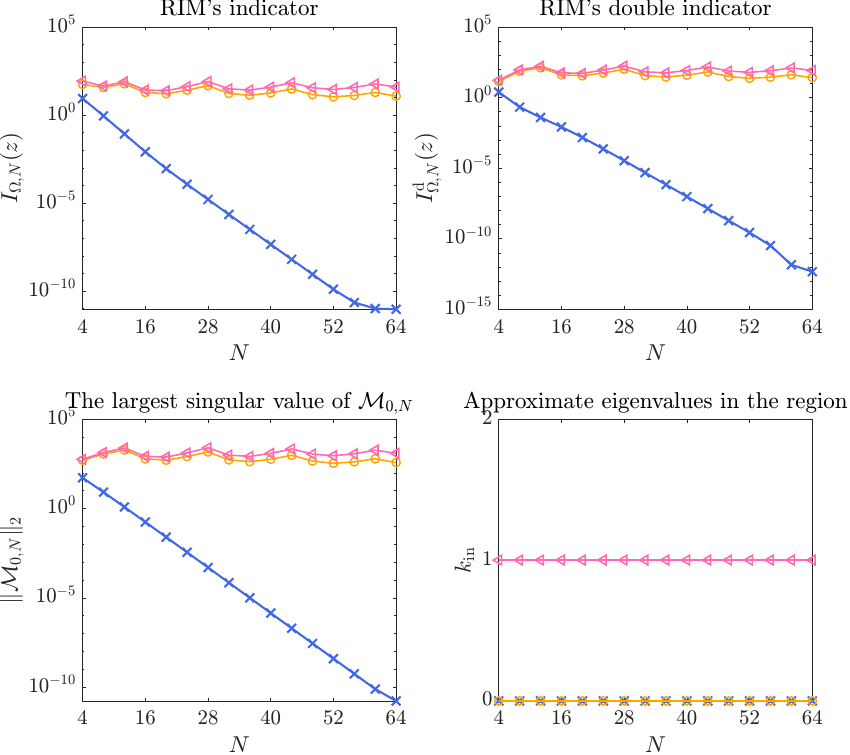}
\includegraphics[width=0.4\linewidth]{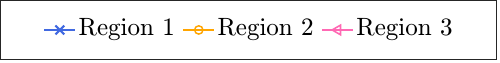}
\caption{The performance of four partitioning criteria on three different
regions in Figure~\ref{fig:3indicator_rg}.
}
\label{fig:3indicator}
\end{figure}

The number of eigenvalues estimated by recursive Beyn's method is
plotted at the bottom-right corner of Figure~\ref{fig:3indicator}.
It can be observed that the proposed algorithm accurately estimates the
number of eigenvalues on all three test regions, even when only a few
quadrature nodes are used.

Three other partitioning methods are also plotted in
Figure~\ref{fig:3indicator} for comparison.
The top-left image is the output of the RIM's indicator \(I_{\Omega,N}(z)\).
To make our experiments more convincing, we also include the RIM's double
indicator
\begin{equation}
\label{eq:rimind2}
I_{\Omega,N}^{\rm{d}}(z)=\left\lVert\frac{1}{2\pi\mi}\sum_{j=0}^{N-1}
\omega_jT(\xi_j)^{-1}\frac{\sum_{j=0}^{N-1}\omega_jT(\xi_j)^{-1}z}
{\Bigl\lVert\sum_{j=0}^{N-1}\omega_jT(\xi_j)^{-1}z\Bigr\rVert_2}\right\rVert_2,
\end{equation}
at the top-right.
This indicator, which is also proposed in~\cite{HSS2016}, performs the
contour integration twice, and is declared to be more robust than
\(I_{\Omega,N}(z)\).
Finally, the largest singular value of \(\mathcal{M}_{0,N}\),
i.e., \(\lVert\mathcal{M}_{0,N}\rVert_2\), is plotted at the bottom-left corner,
which reflects the efficiency of Algorithm~\ref{alg:rim_numrank}.
Remember that Algorithm~\ref{alg:rim_numrank} estimates the number of eigenvalues
in a region by truncating its SVD.
Thus, we expect to see \(\lVert\mathcal{M}_{0,N}\rVert_2\) large on Region
\(3\), and small on Regions \(1\) and \(2\).

Note from Figure~\ref{fig:3indicator} that all three comparative
experiments perform well on Regions \(1\) and~\(3\).
For Region \(1\), all three comparative experiments show an exponential decline to \(0\), indicating that no eigenvalue
is included.
And for Region \(3\), which includes one eigenvalue, the trend in all cases is
non-decaying with respect to the increase of \(N\), which is a signal for
existence of eigenvalues.

Nevertheless, when applied to Region \(2\), these three methods are no longer
reliable.
It can be seen that applying them to Regions \(2\) and \(3\)
gives almost the same result, even though there is no eigenvalue in Region
\(2\).
Comparing on Regions \(1\) and \(2\),
the RIM's indicators and the truncated SVD
may return significantly different values for subregions with the same number of eigenvalues,
which is zero in our example.
Therefore, using a fixed classification threshold \(\dind\) or \(\dnr\) to determine
partitioning is not reliable.

\subsection{Comparison with RSI and RIM}
As comparisons, we use two recently proposed partitioning-based nonlinear
eigensolvers, RIM and RSI, to solve the {\tt{photonics\_1}} problem.

\subsubsection{Using RSI to solve {\tt{photonics\_1}}}
Unlike our algorithm, the reduced subspace iteration
(RSI) algorithm~\cite{TS2024} is a Krylov-based eigensolver.
Instead of directly partitioning the region, it achieves partitioning by
splitting the shift into several different shifts distributed in the
region.
Then, the Arnoldi process is employed on each shift to obtain the
eigenvalues of the whole region.

\begin{figure}[tb!]
\centering
\includegraphics[width=0.7\linewidth]{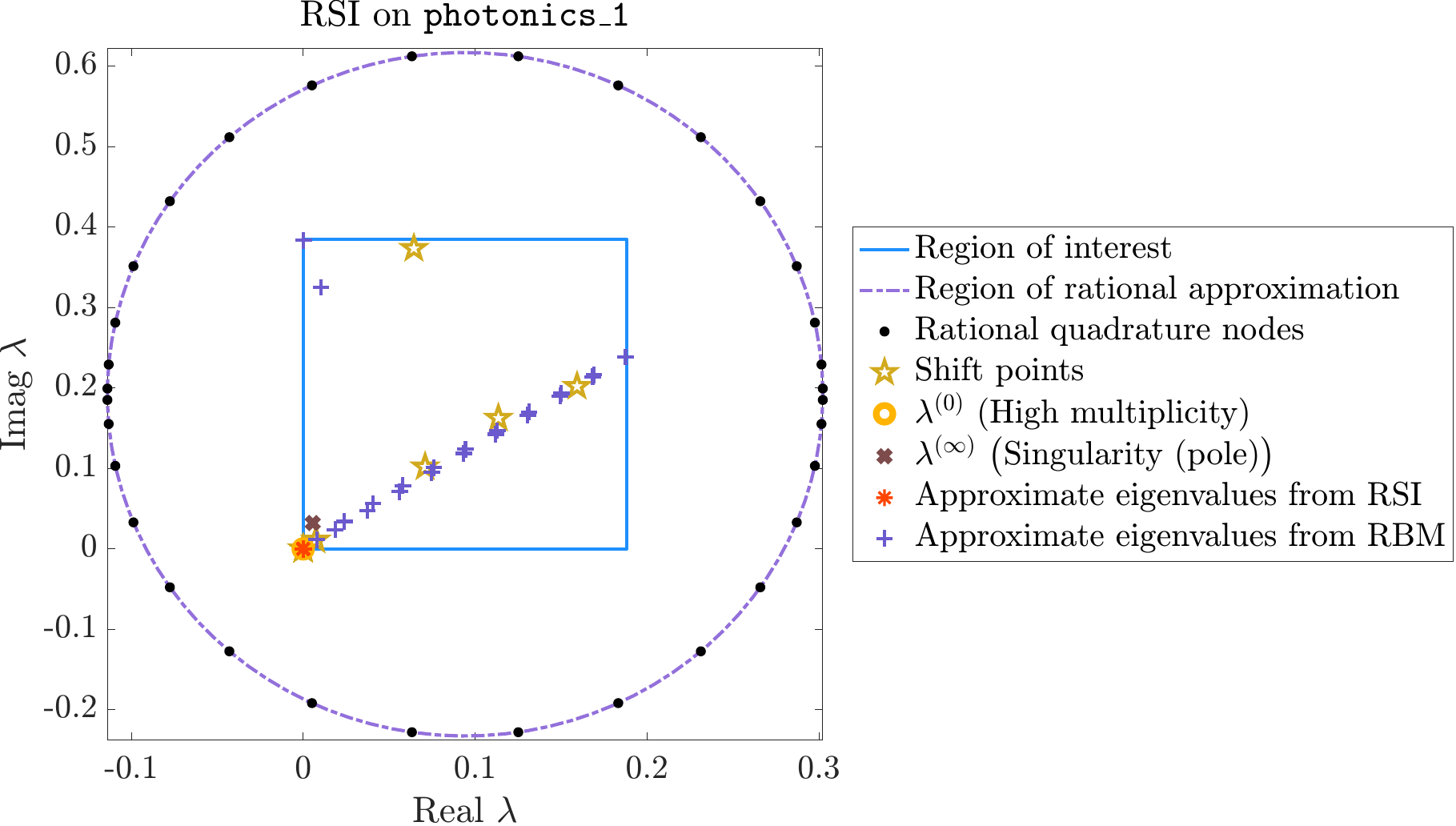}
\caption{Use of the reduced subspace iteration (RSI) algorithm to solve the
{\tt{photonics\_1}} problem.
We used a modified version of the code provided in~\cite{TS2024} to set the
region of rational approximation to be elliptical disks instead of circular
disks.
The convergence threshold is set to \(10^{-12}\).
The region of rational approximation, the rational quadrature nodes, and five shifts
found and used by RSI are marked in the figure.
The approximate eigenvalues generated by recursive Beyn's method (RBM)
are also marked in the figure for comparison.
}
\label{fig:RSI}
\end{figure}

The RSI code is provided in~\cite{TS2024}.
We made slight modifications to the code, enabling the approximate model
of RSI to be applied to elliptical regions.
The result of using RSI to solve the {\tt{photonics\_1}} problem is
illustrated in Figure~\ref{fig:RSI}.%
\footnote{The nonlinear function \(T(\xi)\) is interpolated in RSI.
The region of rational approximation and rational quadrature nodes marked in
Figure~\ref{fig:RSI} are for this use.
For a detailed explanation, see~\cite{TS2024}.}

All \(241\) eigenvalues obtained by RSI are trivial, since they
concentrate around the zero eigenvalue \(\lambda^{(0)}\), which
means that computational resources are severely wasted trying to determine the eigenvalue
with high multiplicity.
Nevertheless, the shifts generated by RSI still capture the
distribution pattern of eigenvalues within the region;
see the shifts and the approximate eigenvalues from RBM in
Figure~\ref{fig:RSI}.
However, the pole \(\lambda^{(\infty)}\) destroys the interpolation
accuracy of RSI, making it fail to calculate the interior eigenvalues
accurately.

\subsubsection{Using RIM to solve {\tt{photonics\_1}} }

\begin{figure}[tb!]
\centering
\includegraphics[width=0.8\linewidth]{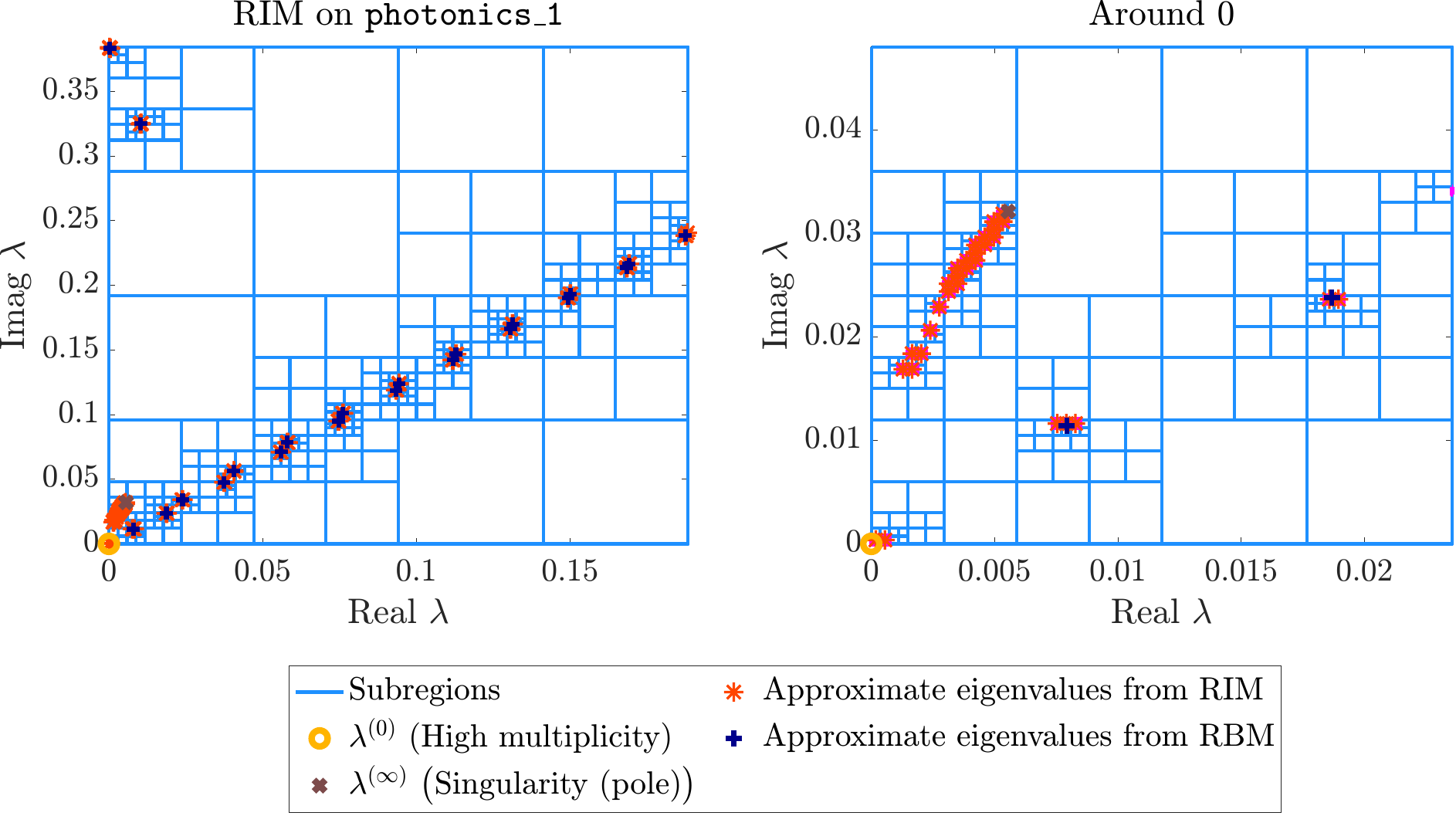}
\caption{Use of the recursive integral method (RIM) to solve the
{\tt{photonics\_1}} problem.
We use RIM’s double indicator~\eqref{eq:rimind2} with a \(32\)-point
Gauss--Legendre quadrature rule.
The convergence threshold is set to \(10^{-3}\), to save computational
resources, while all other parameters match the ones used in~\cite{HSS2016}.
The approximate eigenvalues generated by recursive Beyn's method (RBM)
are also marked in the figure for comparison.
}
\label{fig:RIM}
\end{figure}

The recursive integral method (RIM)~\cite{HSS2016} has already been introduced
in Section~\ref{sub-sec:rim}.
Before presenting the experiments, let us specify the parameters used in our
implementation.
The RIM's double indicator~\eqref{eq:rimind2} is employed with
a Gauss--Legendre quadrature rule on a rectangular contour.
Both the block size \(\ksub=5\) and the number of quadrature nodes on each
segment \(N=32\) are set equal to those used in Figure~\ref{fig:photonics1}
for fairness.
However, the target accuracy is set to \(10^{-3}\) to save computational
resources.
The result is plotted in Figure~\ref{fig:RIM}.

Intuitively, RIM appears to have successfully identified the approximate
locations of all eigenvalues.
In fact, RIM is used for solving similar photonic problems
in~\cite{XS2021}, and shows a satisfactory performance when only a relatively
low precision is needed.
However, a closer look (see Figure~\ref{fig:RIM} right) reveals that two
approximate eigenvalues are calculated around each eigenvalue of
{\tt{photonics\_1}}.
This occurs because the indicator returns large values in every subregion
around an eigenvalue, which means the indicator can hardly locate the
subregions containing eigenvalues accurately.
When comparing to Figure~\ref{fig:photonics1}, it can be observed that the
RIM performs unnecessary partitioning in many regions that do
not contain any eigenvalue, leading to a waste of computational resources.
RIM performs even worse near the pole \(\lambda^{(\infty)}\),
returning many spurious eigenvalues.

\subsection{Implementation with infinite GMRES}
\label{sub-sec:numexp_infgmres}

In this subsection we illustrate the acceleration effect of infGMRES on
recursive Beyn's method.
All three test example, {\tt{gun}}, {\tt{photonics\_1}},
and {\tt{photonics\_2}}, are solved with infGMRES on a uniform grid of
expansion points.
The comparison of time consumed is presented in Table~\ref{tab:time}.

\begin{table}[tb!]
\centering
\caption{Using infinite GMRES to accelerate recursive Beyn's method.
The columns labeled by `Subregions' and `Linear systems' in the table are the
total number of subregions generated and the linear systems solved during the
process of recursive Beyn's method, respectively.
The number of infGMRES expansion points is also shown in the table.
All other parameters, e.g., \(\dmax\), \(N\), and \(\dres\), are the same as
in the experiments in Section~\ref{sub-sec:overall}.
As results, we present the time consumed by infGMRES, the total solving time, and also the number of eigenvalues found.
}
\label{tab:time}
\small
\begin{NiceTabular}{c@{\hspace{0.1em}}c@{\hspace{0.4em}}c@{\hspace{0.1em}}c@{\hspace{0.1em}}c@{\hspace{0.1em}}c@{\hspace{0.1em}}c@{\hspace{0.1em}}c}
\hline
Problem&InfGMRES&Subregions&
\begin{tabular}{c}Linear\\systems\end{tabular}&
\begin{tabular}{c}Expansion\\points\end{tabular}&
\begin{tabular}{c}InfGMRES\\time (s)\end{tabular}&
\begin{tabular}{c}Total\\time (s)\end{tabular}&
\begin{tabular}{c}Eigenvalues\\found\end{tabular}\\
\hline
\multirow{2}*{{\tt{gun}}}
&w/o   & 93  & \phantom{0,0}14,880 & -  & -   & 1066 & 22  \\
&with  & 157 & \phantom{0,0}25,120 & 25 & 132 & 276  & 22\\
\hdottedline
\multirow{2}*{{\tt{photonics\_1}}}
&w/o   & 161 & \phantom{0,0}25,760 & -   & -   & 1488 & 22\\
&with  & 165 & \phantom{0,0}26,400 & 143 & 501 & 1023 & 22\\
\hdottedline
\multirow{2}*{{\tt{photonics\_2}}}
&w/o   & 365  & 3,737,600 & -   & -     & 9340 & 197 \\
&with  & 713  & 7,301,120 & 240 & 2842  & 7646 & \textbf{178}\\
\hline
\end{NiceTabular}
\end{table}

\begin{figure}[tb!]
\centering
\includegraphics[width=0.6\linewidth]{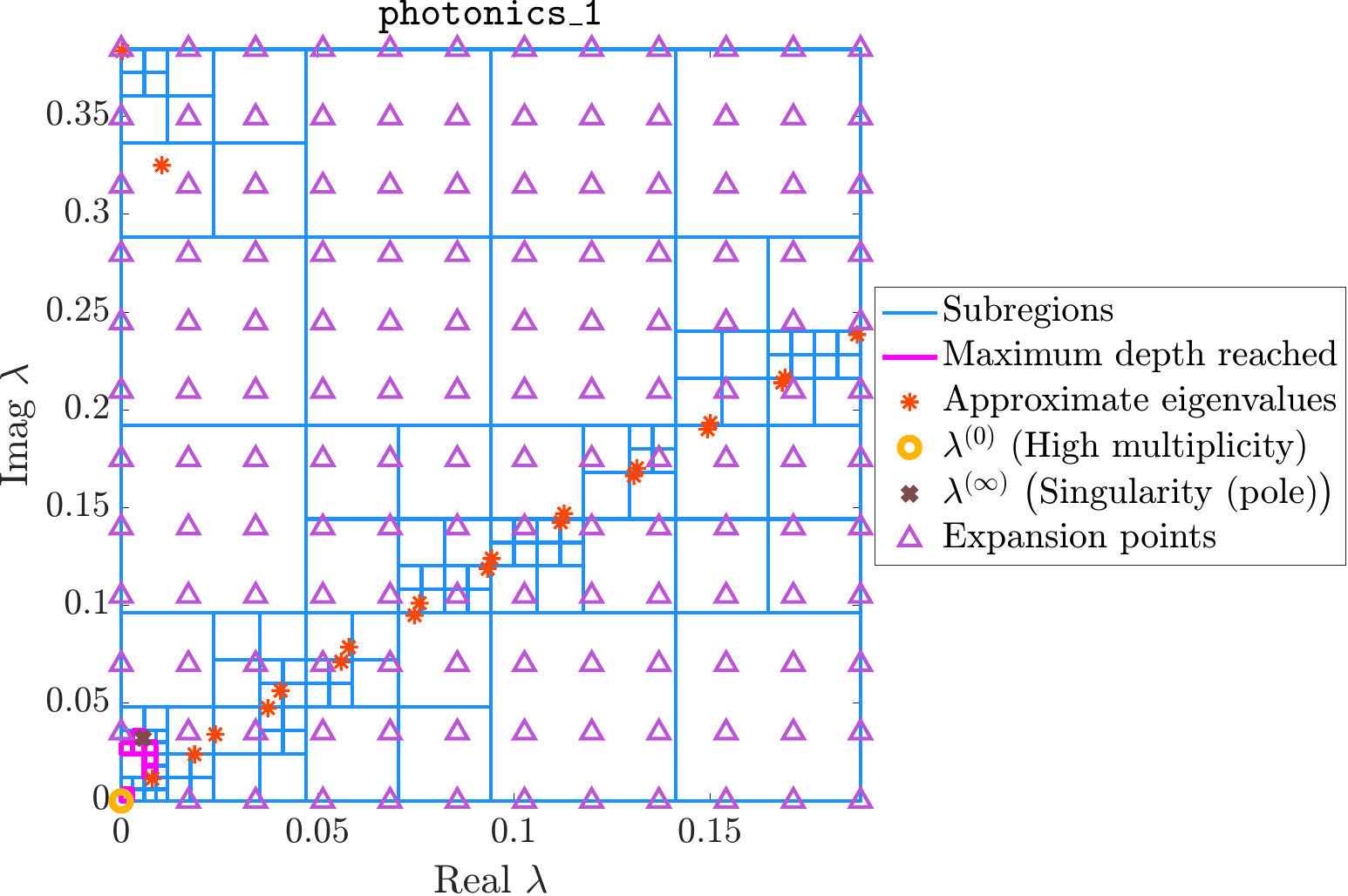}
\caption{Using recursive Beyn's method with infGMRES to solve the
{\tt{photonics\_1}} problem.
The region of interest and the parameters \(\ksub\), \(N\), and
\(\dmax\) are set equal to the ones used in Figure~\ref{fig:photonics1}.
The expansion points of infGMRES are distributed on a \(12\times12\)
uniform grid on the region of interest, except the point at the origin
(\(0\)).
Same as Figure~\ref{fig:photonics1}, all \(22\) eigenvalues are found
and computed with a relative residual below \(10^{12}\).
The difficult points \(\lambda^{(0)}\) and
\(\lambda^{(\infty)}\) are also located.}
\label{fig:photonic_iG}
\end{figure}

As an example, the result of using recursive Beyn's method with infGMRES to solve {\tt{photonics\_1}}
is shown in Figure~\ref{fig:photonic_iG}.
We use \(12\times12\) uniformly distributed expansion points, except
the point at the origin.
Both the approximate eigenvalues and the difficult points
\(\lambda^{(0)}\) and \(\lambda^{(\infty)}\) are found.
However, only \(143\times5=715\) LU decompositions on \(20363\times20363\)
matrices are needed, and the time consumed is only about \(2/3\).

However, as we mentioned before, selecting suitable expansion points to ensure
accurate results from infGMRES is challenging.
As can be observed in Table~\ref{tab:time}, some eigenvalues are \emph{missed} due to the inaccuracy of
infGMRES.
This is partly because infGMRES can hardly solve linear systems near the pole
\(\lambda^{(\infty)}\) accurately.
In fact, it is observed that, even increasing the number of expansion points
helps little to improve the accuracy of the solutions of the linear systems
near singularities.
Therefore, applying infGMRES in intricate regions is technically challenging,
and further investigation is needed to make it robust.

However, since our algorithm continuously subdivides the region
locally, generating a large number of quadrature nodes and linear
systems, there is still significant potential for performance
improvement by infGMRES.
The linear case of this technique is discussed in~\cite{HSY2018}.

\section{Conclusion}
In this work, we propose a contour integral-based, region partitioning
nonlinear eigensolver, i.e., recursive Beyn's method, by combining
Beyn's method and the RIM.
In this algorithm, we replace the RIM's indicator with an estimate of
the number of eigenvalues to determine whether a region should be
further subdivided.
This not only avoids the algorithm's reliance on problem-dependent
parameters, but also improves the accuracy of the region partitioning.
As a result, the algorithm almost never wastefully subdivides regions
that do not contain eigenvalue, which can save computational resources.
Furthermore, unlike the RIM algorithm, which requires a postprocessing stage
to compute eigenvectors, our recursive Beyn's method computes both the
eigenvalues and the corresponding eigenvectors simultaneously.

In terms of experimental results, the proposed algorithm
can accurately identify all eigenpairs, even for the cases when
the region of interest is challenging with singularities,
high-multiplicity eigenvalues, and accumulation points contained,
which can be beneficial for solving certain problems in
quasi-normal mode analysis.

At the end of the paper, we also heuristically introduced the idea of
using infGMRES to accelerate the algorithm.
One of the priorities of our future work will be to adaptively apply
infGMRES to our algorithm, further reducing its computational time
and enhancing its adaptability to various applications.

\section*{Acknowledgments}
Y. Liu and M. Shao were partly supported by the National Natural Science
Foundation of China under grant No.~92370105.
J. E. Roman was partially supported by grants PID2022-139568NB-I00 and RED2022-134176-T funded by MCIN/AEI/10.13039/501100011033 and by ``ERDF A way of making Europe''.
This work was partly carried out while Y.~Liu was visiting Universitat Polit\`ecnica de Val\`encia.

\end{document}